\documentclass[12pt]{article}
\usepackage{amssymb}
\newtheorem{theorem}{Theorem}
\newtheorem{corollary}{Corollary}

\newtheorem{proposition}{Proposition}
\newtheorem{lemma}{Lemma}
\newtheorem{conjecture}{Conjecture}

\title{A Structure Theorem for Positive Density Sets Having the 
Minimal Number of 3-Term Arithmetic Progressions}
\author{Ernie Croot}
\date{\today}

\begin{document}

\maketitle

\begin{abstract}  Assuming the well known conjecture that for any $\gamma > 0$
and $x$ sufficiently large the interval $[x,x+x^\gamma]$ always contains
a prime number, we prove the following unexpected result:  There exist numbers $0 < \rho < 1$ 
arbitrarily close to $0$, and arbitrarily large primes $q$, such that if $S$ is any 
subset of ${\mathbb Z}/q{\mathbb Z}$ of density at least $\rho$, 
having the least number of $3$-term arithmetic progressions among all 
such sets $S$ (of density $\geq \rho$), then there exists an integer
$1 \leq b \leq q-1$ and a real number $0 < d < 1$ (depending only on
$\rho$)  such that
$$
|S\ \cap\ (S+bj)|\ =\ |S|\left (1 - O \left ( {1 \over |\log \rho_1|} \right ) 
\right ),\ {\rm for\ every\ } 0 \leq j < q^d.
$$
This result says that $S$ is ``nearly translation invariant'' in a very 
strong sense.

A curious feature of the proof is that F. A. Behrend's result on large
subsets of $\{1,2,...,x\}$ containing no 3-term arithmetic progressions
is a key ingredient.  The proof also uses exponential sums, a result of K. F.
Roth, as well as a result of P. Varnavides. 

\end{abstract}

\section{Introduction}

Given a subset $S$ of ${\mathbb Z}/q{\mathbb Z}$, we will say that $S$
has density $\rho$ if $|S| = \rho q$.  A well known result of K. F. Roth 
\cite{roth} asserts that for $q$ sufficiently
large, and for $(\log\log q)^{-1} < \rho \leq 1$, any subset of 
${\mathbb Z}/q{\mathbb Z}$ having density at least $\rho$, must contain a
non-trivial three term arithmetic progression; that is, a triple of numbers,
$a,b,c$, $ a \not \equiv b \pmod{q}$, satisfying 
\begin{equation} \label{abc}
a + b\ \equiv\ 2c \pmod{q}.
\end{equation}
E. Szemeredi \cite{szemeredi}, R. Heath-Brown \cite{heathbrown}, and J. Bourgain
\cite{bourgain} have improved considerably on Roth's result, with the Bourgain's
work being the most recent:
\bigskip

\begin{theorem} {\bf (J. Bourgain)} \label{bourgain_result} 
For $q$ sufficiently large and 
$$
B \sqrt{\log\log q \over \log q}\ <\ \rho\ \leq\ 1,\ {\rm where\ }
B > 0\ {\rm is\ some\ constant},
$$ 
any subset of ${\mathbb Z}/q{\mathbb Z}$ of density $\rho$
contains a triple $a,b,c$ (with $a \neq b \neq c$) satisfying (\ref{abc}).
\end{theorem}
\bigskip

As a consequence of Roth's result mentioned above, Varnavides \cite{varnavides}
proved the following
\bigskip

\begin{theorem} {\bf (P. Varnavides)} \label{varnavides_theorem}
Given $0 < \rho \leq 1$, and a sufficiently large prime $q$, 
if $S$ is any subset of 
${\mathbb Z}/q{\mathbb Z}$ of density $\rho$, there exists $\kappa > 0$
so that
$$
\mu_q(S)\ =\ {\#\{a,b,c \in S\ :\ a+b \equiv 2c \pmod{q}\} \over q^2}\ >\ 
\kappa.
$$
In fact, one can take
$$ 
\kappa\ =\ {\rho \over 16 h(\rho/2)^2},
$$
where $h(\rho)$ is defined to be the least integer such that if 
$m \geq h(\rho)$ and $T$ is any subset of $\{1,...,m\}$ with
density at least $\rho$, then $T$ contains a three term arithmetic progression.
\end{theorem}
The version of the theorem appearing in \cite{varnavides} does not
give such an explicit value for $\kappa$, and the result is stated in
terms of subsets of $\{1,2,...,x\}$, rather than ${\mathbb Z}/q{\mathbb Z}$.
For these reasons, we give our own proof of this result in section
\ref{varnavides_section}.

Although the function $\mu_q(S)$ counts both non-trivial and trivial
solutions to $a + b \equiv 2c \pmod{q}$ (trivial means $a \equiv b 
\equiv c \pmod{q}$), if $q$ is large enough in terms of $\rho$, we can
deduce that
\begin{equation} \label{nontrivial}
{\#\{a,b,c \in S,\ a \not \equiv b \pmod{q}\ :\ a + b \equiv 2c \pmod{q}\}
\over q^2}\ >\ \kappa, 
\end{equation}
since there can be at most $|S| < q$ triples $a,b,c\in S$ 
with $a \equiv b \equiv c \pmod{q}$; and so, the contribution of these
trivial solutions to $\mu_q(S)$ is thus only $O(1/q)$, which tends to $0$
as $q$ tends to infinity, which thus proves (\ref{nontrivial}).

Combining Varnavides's and Bourgain's results, one can show that if 
$S$ has density $\geq \rho$ modulo $q$, and if $q$ is a sufficiently large prime,
then

\begin{equation} \label{bourgain_varnavides}
\mu_q(S)\ \geq\ \exp \left ( -D {|\log \rho| \over \rho^2} \right )
\end{equation}
where $D >0$ is an absolute constant.
\bigskip

It is of interest to try to find, for a given density $0 < \rho \leq 1$, the 
smallest value of $\kappa$ so that the conclusion of Varnavides's theorem 
still holds.  This smallest $\kappa$ is given by the following 
function:
$$
r_q(\rho)\ =\ \min_{S \subseteq {\mathbb Z}/q{\mathbb Z} \atop |S| \geq \rho q}
\mu_q(S).
$$
One approach to understanding $r_q(\rho)$ is to try to understand 
the structure of {\it critical sets}, which are subsets $S$
of ${\mathbb Z}/q{\mathbb Z}$ of density $\rho$ and 
$$
\mu_q(S)\ =\ r_q(\rho).
$$

It would seem that the problem of understanding the structure of these 
critical sets is a much more difficult problem than that of understanding  
the behavior of $r_q(\rho)$; and so, it would seem that we should
try to make progress on the order of growth of $r_q(\rho)$ by some 
other method {\it first} before trying to tackle questions about such sets.
Even so, one might would think that they have very little
structure; however, motivated by the main result of this 
paper (theorem \ref{main_theorem}) we have the following
conjecture, which claims, on the contrary, that these sets have 
a considerable amount of additive structure:
\bigskip

\begin{conjecture}  Given $0 < \rho < 1$, and $q$ sufficiently large,
there exists $0 < d < 1$ such that
if $S$ is a critical set of ${\mathbb Z}/q{\mathbb Z}$ of density $\rho$,
then there exists a number $1 \leq b \leq q-1$ such that 
$$
|S\ \cap\ (S+jb)|\ >\ (1 - n(\rho)) |S|,\ {\rm for\ every\ }
0 \leq j < q^d,
$$
where $n(\rho)$ is a function of $\rho$ only, and which tends to $0$ as 
$\rho$ tends to $0$.
\end{conjecture}
\bigskip

The main theorem of the paper is a weakened version of this conjecture.
Before we state it, first define
$$
r(\rho)\ =\ \liminf_{q \geq 3\atop q\ {\rm prime}} r_q(\rho).
$$
We will need the following, famous conjecture from prime number theory:
\bigskip

\begin{conjecture} \label{primes}  Given $\theta > 0$, we have for all $x$ 
sufficiently large that $[x,x+x^\theta]$ contains a prime number.  
\end{conjecture}
From \cite{baker} the conjecture is known to hold for all $\theta > 0.525$.
\bigskip

Our main theorem is as follows:
\bigskip

\begin{theorem} \label{main_theorem}  Assume conjecture \ref{primes} holds. 
For any $0 < \rho_0 < 1$, there exist numbers
$0 < \rho_1 < \rho_0$ and $0 < d < 1$, and infinitely many primes $q$
such that the following holds:  If 
$S \subseteq {\mathbb Z}/q{\mathbb Z}$ has density $\geq \rho_1$,
and has the least number of $3$-term arithemtic progressions modulo $q$
among all sets of density $\geq \rho_1$ (that is, $\mu_q(S) = r_q(\rho_1)$),
then there exists a number $1 \leq b \leq q-1$ such that 
$$
|S \cap (S+bj)|\ \geq\ |S| \left ( 1 - O \left ( {1 \over |\log \rho_1|} 
\right ) \right ),\ {\rm for\ every\ } 0 \leq j \leq q^d,
$$
where the implied constant in the big-O is absolute.
\end{theorem}
\bigskip

If we could prove the following conjecture on the behavior of 
$r_q(\rho)$, then we could remove the assumption that Conjecture
\ref{primes} holds in the above theorem:\footnote{At the end of the proof
of Theorem \ref{main_theorem}, in section \ref{main_theorem_section},
we show how the truth of this conjecture implies this stronger
version of our main theorem}
\bigskip

\begin{conjecture} \label{rqn_conjecture}  
There exists a constant $A > 1$ such that 
for every $0 < \rho \leq 1$, every $0 < \theta < 1$, and
$n$ sufficiently large, 
$$
{r_q(\rho) \over r_n(\rho)}\ \in\ \left [ {1 \over A},\ A \right ],
\ {\rm for\ every\ }q \in [n,n+n^\theta].
$$
\end{conjecture}
\bigskip

We now briefly mention the main ideas and ingredients in the proof of 
Theorem \ref{main_theorem}, by describing how to prove a certain toy version 
of the theorem:  For a given set of integers $T$ having
density at least $\rho$ modulo $q$, define the exponential sum
$$
f_T(t)\ =\ \sum_{s \in T} e(st),\ 
$$
where $e(u) = \exp(2\pi i u)$.  Using the identity
$$
{1 \over q} \sum_{j=0}^{q-1} e\left ( {ja \over q} \right )\ =\ 
\left \{ \begin{array}{rl} 
1,\ & {\rm if\ }q\ {\rm divides\ } a;  \\
0,\ & {\rm if\ }q\ {\rm does\ not\ divide\ } a,
\end{array} \right.
$$
one can easily show that 
\begin{eqnarray} \label{T_sum}
q^2 \mu_q(T)\ &=& \#\{a,b,c \in T\ :\ a+b \equiv 2c \pmod{q}\}\nonumber \\
&=&\ {1 \over q} \sum_{|a| < q/2} f_T \left ( {a \over q} \right )^2 
f_T \left ( {-2a \over q} \right ).
\end{eqnarray}
It turns out that to estimate the sum on the right to within an error of size 
less than $\epsilon |T| q < \epsilon q^2$, we need only sum over values $a$ 
in a set with at most $(\rho \epsilon^2)^{-1}$ elements.  Now, given
a critical set $S$ of density $\rho$, we show that we can multiply
the set by an integer (not divisible by $q$), to produce a new 
set $W$ modulo $q$, where we think of $W$ as being a subset of 
$\{0,1,...,q-1\}$, such that if
$$
\left | f_W \left ( {a \over q} \right ) \right |\ \geq\ \epsilon |W|,
\ \left | {a \over q} \right | < {1 \over 2},
$$
then 
\begin{equation} \label{aq_assumption}
\left | {a \over q} \right | < q^{- \rho\epsilon^2},\ {\rm and\ }
a\ {\rm is\ divisible\ by\ }4.
\end{equation}

Now suppose that $p$ is a prime number very close to $q/2$.  
As a consequence of (\ref{aq_assumption}), we show that the sum on the 
right-hand-side of (\ref{T_sum}) is, for small values of $\epsilon$, 
\begin{eqnarray} \label{pq_equation}
q^2 \mu_q(T)\ &\approx& {1 \over q} \sum_{|a| < p/2} f_W \left ( {a \over p} \right )^2
f_W \left ({-2a \over p} \right )\nonumber \\
&=&\ {p \over q}\ \#\{a,b,c \in W\ :\ a+b \equiv 2c \pmod{p}\},
\end{eqnarray}
where, in this context, we will use the notation $u(x) \approx v(x)$
to loosely mean that the functions $u(x)$ is ``extremely close'' to 
$v(x)$; and, as $\epsilon \to 0$ and $q \to \infty$, the ratio of 
these two functions tends to $1$.

Next, we show that for certain special values of $\rho$, 
the proportion of residue classes modulo $p$ 
occupied by $W$ is very nearly $\rho$; so, most residue classes
modulo $p$ contain either $0$ or $2$ elements of $W$.  
A consequence of this is that if $V_0$ is the set of all
elements of $W$ that are in $[0,p-1]$, and if $V_1$ is the set of
all elements of $W$ in $[p,q-1]$, then 
$$
\left |V_0\ \cap\ (V_1 - p) \right |\ \geq\ |V_0| \left ( 1 -
O \left ( {1 \over |\log \rho|} \right ) \right ).
$$ 
Using an additional trick that involves rotating the set $W$ modulo $q$ (translating
by an integer $k$), we can show that 
$$
| W\ \cap\ (W-p) |\ \geq\ |W| \left ( 1 - O \left ( {1 \over |\log \rho|} 
\right ) \right ),
$$
and so there exists an integer $b$, not divisible by $q$, so that 
$$
| S\ \cap\ (S+b) |\ \geq\ |S| \left ( 1 - O \left ( {1 \over |\log \rho|}
\right ) \right ),
$$
since the elements of $S$ are multiples of elements of $W$ modulo $q$.  
The additional parameter $j$ in the intersection $S\ \cap\ (S+jb)$ stated
in theorem \ref{main_theorem} comes about via only a slight generalization
of the above argument. 

The way we deduce that the proportion of resiude classes modulo $p$ occupied
by $W$ is nearly the same as the proportion modulo $q$ is as follows.  First, 
we need the following result of F. A. Behrend:
\bigskip

\begin{theorem} \label{behrend_theorem}  For sufficiently large $q$ there exists
a subset $S$ of ${\mathbb Z}/q{\mathbb Z}$ with
$$
|S|\ >\ {q \over \exp(C \sqrt{\log q})},\ {\rm for\ some\ constant\ }C > 0,
$$
such that $S$ contains no three term arithmetic progressions modulo $q$.
\end{theorem}
\bigskip

\noindent As a consequence of this theorem, we show that

\begin{corollary} \label{behrend_consequence}
Given $0 < \rho < 1$, we have that 
$$
r(\rho)\ \leq\ \left ( -{1 \over C^2} \log^2 (4\rho) \right ).
$$
\end{corollary}
So, Behrend's result shows that the function $r(\rho)$ decays quite
rapidly compared with $\rho$.
\bigskip

Now, $W$ must occupy at least $\rho$ of the residue 
classes modulo $p$.  As a consequence of the above corollary, one 
can show that for certain values of $\rho$, if $W$ had density
more than
$$
\rho \left ( 1 + {G \over |\log \rho|}\right ),
$$
modulo $p$, where $G > 0$ is some particular constant,
then the number of solutions $a+b \equiv 2c \pmod{p}$,
$a,b,c \in W$, would have to be a large multiple of the number of
such solutions modulo $q$, and so we would have that (\ref{pq_equation})
could not hold.
\bigskip

The rest of the paper is orgainized as follows.  In the next section,
we give a proof of theorem \ref{main_theorem}, as well as an indication
of how conjecture \ref{rqn_conjecture} can be used to replace 
the assumption that conjecture \ref{primes} holds.  In several of 
the sections after the proof of the main theorem, we prove several
propositions\footnote{In my terminology, ``proposition'' is synonymous with
``large lemma'' or perhaps ``meta-theorem''; that is, a proposition is
a result whose proof is either too large, or is too technical to be considered
a lemma, yet is not sufficiently general or interesting to be considered
a theorem.} that are 
 used in the proof of theorem \ref{main_theorem},
as well as in the proofs of other propositions.  Finally, in section
\ref{lemma_section}, we state and prove several technical lemmas
and corollaries that appear throughout the paper. 
\bigskip

\section{Proof of Theorem \ref{main_theorem}.} \label{main_theorem_section}
\bigskip

We first pin down the value of $\rho_1$ for which we will prove that 
$S$ satisfies the conclusion of the theorem, and the following
proposition gives us the answer we seek:

\begin{proposition} \label{rho_prop}
Suppose that $0 < \rho_0 < 1$ and $k>1$.  Then, there exists a number
$0 < \rho_1 < \rho_0$ such that for all 
$$
\rho\ >\ \rho_1 \left ( 1 + {2C^2\log k \over |\log \rho_1|} \right ),
$$
we will have
$$
r(\rho)\ >\ k\ r(\rho_1),
$$
where $C$ is the same constant that appears in 
Theorem \ref{behrend_theorem} mentioned in the introduction.
\end{proposition}
The proof of this proposition appears in section \ref{rho_prop_proof}.

We will assume henceforth that $\rho_1$ is any number satisfying the
conclusion of this proposition for $k=100$.  Then, for 
$0 < \epsilon < 1$, to be chosen later, let $q$ be any 
sufficiently large prime so that 
\bigskip

\begin{equation} \label{q_assumption}
r(\rho_1)\ >\ {r_q(\rho_1) \over 2},
\end{equation}
and so that for every prime 
$$
p\ \in\ \left [ { q \over 2},\ {q \over 2} + q^{\rho_1 \epsilon^2/2} \right ]
$$ 
we have 
\begin{equation} \label{T_assumption}
r_p\left ( \rho_1 \left ( 1 + {2C^2 \log 100 \over |\log \rho_1|} \right )
\right )\ >\ {1 \over 2}\ r \left ( \rho_1 \left ( 1 + {2C^2 \log 100 \over |\log \rho_1|}
\right ) \right ).
\end{equation}
Note:  Since we have assumed conjecture \ref{primes}, we have
that the interval above always contains a prime for sufficiently large $q$.  
We also note that every sufficiently large prime $p$ will satisfy 
(\ref{T_assumption}), and infinitely
many primes $q$ will satisfy (\ref{q_assumption}).  Thus, we will have that
both (\ref{q_assumption}) and (\ref{T_assumption}) hold for infinitely many 
primes $q$.
\bigskip

Define the exponential sum
$$
f_S(t)\ =\ \sum_{s \in S} e(st),
$$
where $e(u) = \exp(2\pi i u)$.  We first claim that
$$
\#\left \{ a \in {\mathbb Z},\ |a| < {q \over 2}\ :\ 
\left |f_S\left ( {a \over q} \right )\right | > 
\epsilon |S| \right \}\ \leq\ {1 \over \rho_1 \epsilon^2}.
$$
To see this, we have from Parseval's identity (lemma 
\ref{parseval} in section \ref{lemma_section}) that 
$$
\sum_{|a|<q/2} \left |f_S\left ( {a \over q} \right )\right |^2\ =\ q |S|.
$$
Now, if there were more than $(\rho_1 \epsilon^2)^{-1}$ values of 
$a$ for which $|f_S(a/q)| > \epsilon |S|$, then the sum on the left hand
side would exceed $|S|^2/\rho = q |S|$, and so couldn't equal $q |S|$
as claimed. 

We next need the following lemma:

\begin{lemma} \label{cube_lemma}
There exists an integer $1 \leq h \leq q-1$ such that the set
$$
S'\ =\ h S\ =\ \{h s \pmod{q}\ :\ s \in S\} 
$$ 
has the property that
\begin{equation} \label{S'_property}
\left |f_{S'}\left ( {a \over q} \right )\right | > \epsilon |S|,\ 
|a| \leq q/2\ \Longrightarrow\ |a| < q^{1 - \rho_1 \epsilon^2}.
\end{equation}
\end{lemma}

The proof of this lemma makes use of the pigeonhole principle and 
can be found in section \ref{lemma_section}.  Now, for integers 
$k$ and $1 \leq v < q^{\rho_1 \epsilon^2/4}$, to be chosen later, let 
\begin{eqnarray} \label{W_definition}
W\ &=&\ W(k,v)\ =\ (4v)^{-1}S' + k\ =\ 
\{ (4v)^{-1} s + k \pmod{q}\ :\ s \in S'\}\nonumber \\
&=&\ \{ (4v)^{-1} h s + k \pmod{q}\ :\ s \in S\}. 
\end{eqnarray}
We note that the number of solutions to $a+b \equiv 2c \pmod{q}$, $a,b,c \in W$,
is the same as the number of solutions among the elements of $S$; also, one 
can easily see that for $|a| \leq q/2$,  
\begin{equation} \label{S''_equation}
\left | f_{S'} \left ( {(4v)^{-1} a \over q} \right ) \right |\ =\ 
\left | f_W \left ( {a \over q} \right ) \right |
> \epsilon |S|\ \Longrightarrow\  |a| < 4vq^{1-\rho_1 \epsilon^2},\ 
{\rm and\ } 4v|a.
\end{equation}

We now require the following proposition:

\begin{proposition} \label{lambda_prop}
Suppose $p \in [q/2,q/2+q^{\rho_1 \epsilon^2/2}]$ and that
$1 \leq v \leq q^{\rho_1 \epsilon^2/4}$ is an integer.  
Then, for $W = W(k,v)$, we have 
$$
\mu_q(W) \ =\ {p^3 \over q^3} \mu_p(W)
\ +\ 
O \left ( \left (  {\epsilon \over \rho_1} \right )^{1/3} \rho_1^2 \right ), 
$$
where the constant in the big-O is absolute.
\end{proposition}
The proof of this result can be found in section \ref{lambda_section}.

Now suppose that $R(k,v)$ is the set of residue classes modulo $p$ that are 
occupied by $W(k,v)$.  If
$$
|R(k,v)|\ >\ \rho_1 \left ( 1 + {2C^2 \log 100 \over |\log \rho_1|} \right ) p,
$$
then we would have from (\ref{q_assumption}), (\ref{T_assumption}), 
and proposition \ref{lambda_prop} that 
\begin{eqnarray}
r(\rho_1)\ &>&\ {r_q(\rho_1) \over 2}\ =\ {\mu_q(W) \over 2} \nonumber \\
&=&\ {p^3 \mu_p(W) \over 2q^3}  
\ +\ O \left ( \left ( {\epsilon \over \rho_1} \right )^{1/3} \rho_1^2 \right ) 
\nonumber \\
&\geq&\ {p^3 \over 2q^3} 
r_p\left ( \rho_1 \left ( 1 + {2C^2 \log 100 \over |\log \rho_1|} \right )
\right )\ +\ O \left ( \left ({\epsilon \over \rho_1} \right )^{1/3} \rho_1^2 \right )
\nonumber \\
&>&\ {p^3 \over 4q^3}
r \left ( \rho_1 \left ( 1 + {2C^2 \log 100 \over |\log \rho_1|} \right )
\right )\ +\ O \left ( \left ( {\epsilon \over \rho_1} \right )^{1/3} 
\rho_1^2 \right ) \nonumber \\
&>&\ {100 p^3 \over 4q^3} r(\rho_1)\ +\ 
O \left ( \left ( {\epsilon \over \rho_1} \right )^{1/3} \rho_1^2 \right ),
\nonumber
\end{eqnarray}
which is impossible once $\epsilon > 0$ is small enough.  Thus, we conclude
\begin{eqnarray} 
{|S| \over 2}\ \leq\ |R(k,v)|\ &\leq&\ \rho_1 \left ( 1 + {2C^2 \log 100 \over 
|\log \rho_1|} \right ) p \nonumber \\
&=&\ |S|\left ({1 \over 2}\ +\ O \left ( {1 \over |\log \rho_1|} \right )
\right ); \nonumber
\end{eqnarray}
and so,
\begin{equation} \label{R_upper_bound}
|R(k,v)|\ =\ |S| \left ( {1 \over 2} + O \left ( {1 \over |\log \rho_1|} \right )
\right ).
\end{equation}
This tells us that a typical residue class modulo $p$ in $R(k,v)$
contains either $0$ or $2$ elements of $W(k,v)$; more precisely, it says that
the number of progressions modulo $p$ containing only one element of 
$W(k,v)$ is $O(|S|/|\log \rho_1|)$.  Thus, if we let 
$W_0(k,v)$ denote the integers in $W(k,v)$ that lie in $[0,p-1]$, and
if we let $W_1(k,v)$ denote the integers in $W(k,v)$ that lie in 
$[p,q-1]$, then our equation for $R(k,v)$ says that
\begin{equation} \label{W0W1_equation}
|W_0(k,v)\ \cap\ (W_1(k,v)-p)|\ =\ |S| \left ( {1 \over 2} 
\ -\ O \left ( {1 \over |\log \rho_1|} \right ) \right ).
\end{equation}

We will now show that
\begin{equation} \label{W0_equation}
|W(0,v)\ \cap\ W(-p,v)|\ \geq\ |S| \left ( 1\ -\ O \left ( {1 \over 
|\log \rho_1|} \right ) \right ).
\end{equation}
Here is the proof: First, suppose that $w \in W_0(0,v)\ \cap\ (W_1(0,v)-p)$, 
and note that $0 \leq w \leq p-1$.  For such $w$ we will have
$w \in W(0,v)\ \cap\ W(-p,v)$.  Next, suppose that 
$w - p \in W_0(-p,v)\ \cap\ (W_1(-p,v)-p)$, and note that 
$w \geq p$.  For these integers $w$ we will also have that 
$w \in W(0,v)\ \cap\ W(-p,v)$.  Since the two sets of integers $w$ considered
are disjoint, we must have
\begin{eqnarray}
|W(0,v)\ \cap\ W(-p,v)|\ &\geq&\ |W_0(0,v)\ \cap\ (W_1(0,v)-p)|\ \nonumber \\
&&\ \ \ \ \ \ \ +\ |W_0(-p,v)\ \cap\ (W_1(-p,v)-p)|, \nonumber
\end{eqnarray}
and so (\ref{W0_equation}) follows from this inequality and 
(\ref{W0W1_equation}). 

From this it follows that
\begin{eqnarray}
|S\ \cap\ (S - (4v)h^{-1} p)|\ &=&\ |W(0,v)\ \cap\ W(-p,v)| \nonumber \\
&\geq&\ |S| \left ( 1 - O \left ( {1 \over |\log \rho_1|} \right )\right ), 
\nonumber
\end{eqnarray}
which proves the theorem, since we can take $v$ to be any integer
in $[1,q^{\rho_1 \epsilon^2/4}]$.  Note that here we are thinking of the set
$S - (4v)h^{-1} p$ as a set of integers in $\{0,1,...,q-1\}$.
\bigskip

The reason that conjecture \ref{rqn_conjecture} allows us to remove 
the assumption that conjecture \ref{primes} holds is as follows:  
The proof of the above theorem actually shows more than is stated.
It shows that if $S$ satisfies 
the hypotheses of the theorem, and if $q$ is a prime satisfying
$$
r(\rho_1)\ >\ {r_q(\rho_1) \over 2},
$$
and if $n \in (q/2,q/2 + q^{\rho_1 \epsilon^2/2})$ is some 
integer all of whose non-trivial divisors are greater than 
$q^{\rho_1 \epsilon^2/4}$, and which satisfies
\begin{equation} \label{r_equation}
r_n(\rho_2)\ >\ F\ r(\rho_2),\ {\rm where\ } 0 < F \leq 1
\ {\rm does\ not\ depend\ on\ }\rho_2,\
\end{equation}
where
$$
\rho_2\ =\ \rho_1 \left ( 1 + {2C^2 \log 100 \over |\log \rho_1|} \right ),
$$
then the conclusion of the theorem holds if $q$ is sufficiently large.
We know that this last condition (\ref{r_equation}) 
holds for all sufficiently large
{\it primes} $n$, but it is not clear that it holds for all
sufficiently large integers.

Now, if we assume conjecture \ref{rqn_conjecture}, then for every
$0 < \rho < 1$ we will have that there exists $A \geq 1$ so that
\begin{equation} \label{rn_eqn}
r(\rho)\ =\ \liminf_{q \geq 3 \atop q\ {\rm prime}} r_q(\rho)
\ \leq\ A\ \liminf_{n \geq 3} r_n(\rho).
\end{equation}
To see this, we first note from the result in \cite{baker}, the 
interval $[n,n + n^{0.53}]$ contains a prime number.  For each
such $n$, let $p_n$ denote a prime in this interval.  Then, from
conjecture \ref{rqn_conjecture}, we get that there exists an $A > 0$
(which does not depend on $n$ or $\rho_1$), such that
$$
r_{p_n}(\rho)\ \leq\ A\ r_n(\rho),
$$
and this implies the inequality (\ref{rn_eqn}) above.

From (\ref{rn_eqn}) one now sees that (\ref{r_equation}) holds for all 
sufficiently large integers $n$ (prime or not); and so, our assertion
that conjecture \ref{rqn_conjecture} can be used to remove the 
assumption conjecture \ref{primes} follows.
\bigskip

\section{Proof of Theorem \ref{varnavides_theorem}} \label{varnavides_section}
\bigskip

Suppose $q > h(\rho/2)$.
Set $k = \lfloor h(\rho/2)\rfloor + 1$, and let $P$ be the set of all 
$k$-term arithmetic progressions $a,a+d,a+2d,...,a+(k-1)d$ modulo $q$,
$d$ is not $0$ modulo $q$.  We treat the progression $a+(k-1)d, ..., a$ as
distinct from $a,...,a+(k-1)d$.  Clearly, 
$$
|P|\ =\ q(q-1),
$$
since there are $q$ choices for $a$ and $q-1$ choices for $d$.

Given a progression $h \in P$, we let
$$
R(h)\ =\ {|h\ \cap\ S| \over k}.
$$
If $R(h) \geq \rho/2$, then $h \cap S$ must contain a non-trivial 
$3$-term arithmetic progression, from the way we have defined $h(\rho)$ 
(and $h(\rho/2)$).  When $q$ is sufficiently large, such a progression can
obviously only be a subset of at most $k^2$ progressions $h \in P$; and so,
for $q$ sufficiently large,
\begin{eqnarray} \label{mu_equation}
q^2 \mu_q(S)\ &=&\ \#\{a,b,c \in S\ :\ a+b \equiv 2c \pmod{q}\}
\nonumber \\
&\geq&\ {\#\{h \in P\ :\ R(h) \geq \rho/2\} \over k^2}.  
\end{eqnarray}
To bound this last quantity from below, we first note that 
\begin{eqnarray}
\sum_{h \in P} R(h)\ &=&\ {1 \over k} \sum_{h \in P} |S\ \cap\ h|\ =\ 
{1 \over k} \sum_{s \in S} \sum_{h \in P \atop s \in h} 1 \nonumber \\ 
&=&\ {1 \over k} \sum_{S \in S} k (q-1) \nonumber \\
&=&\ |S| (q-1). \nonumber
\end{eqnarray} 
The second to the last line follows since each $s \in S$ can be in any one
of the $k$ terms of a $k$-term progression (hence the factor $k$); and,
once this term is specified, there are $q-1$ choices for the common 
difference of such a progression (that $s$ lies in). 

Now, we get
\begin{eqnarray}
\sum_{h \in P \atop R(h) \geq \rho/2} 1\ &\geq&\ 
\sum_{h \in P} \left ( R(h) \ -\ {\rho \over 2} \right ) \nonumber \\
&\geq&\ |S| (q-1)\ -\ {\rho \over 2} |P| \nonumber \\
&=&\ \rho q(q-1) - {\rho \over 2} q(q-1) \nonumber \\
&\geq&\ {\rho q^2 \over 4}. \nonumber
\end{eqnarray}
Combining this with (\ref{mu_equation}) now gives
$$
\mu_q(S)\ \geq\ {\rho \over 4 k^2}\ \geq\ {\rho \over 16 r(\rho/2)^2},
$$ 
which proves the theorem.
\bigskip

\section{Proof of Proposition \ref{rho_prop}} \label{rho_prop_proof}
\bigskip

As a consequence of theorem \ref{behrend_theorem} and 
corollary \ref{behrend_consequence}, both of which appear in
the introduction, we have the following lemma, which we will need
for our proof:

\begin{lemma} \label{further_consequence}
For any sequence of numbers $x_1, x_2, x_3, ...$ in $[0,1]$ that converges
to $0$, there are infinitely many integers $n$ such that
$$
{r(x_{n+1}) \over r(x_n)}\ <\ \exp \left (-{1 \over 2C^2} (\log^2 x_{n+1} - 
\log^2 x_n ) \right ).
$$ 
\end{lemma}

Given $k \geq 1$ define the sequence
$$
x_1 = \exp(-2 C^2 \log k),\ \ 
x_{n+1} = x_n \left ( 1 - {C^2\log k \over |\log x_n|} \right ).
$$
We note that $x_1$ was chosen so that the following holds for all
$n \geq 1$:
$$
x_n\ <\ x_{n+1} \left ( 1 + {2 C^2 \log k \over |\log x_{n+1}|} \right ),
$$
for $x_n$ sufficiently close to $0$.

Clearly this sequence tends to $0$ and lies in $[0,1]$.  Thus,
the conditions of the above lemma are satisfied; and so
there exist terms $x_n$ arbitrarily close to $0$ such that
\begin{eqnarray}
{r(x_{n+1}) \over r(x_n)}\ &<&\ 
\exp \left (-{1 \over 2C^2} \left ( \log^2 \left ( x_n 
\left ( 1 - {C^2\log k \over |\log x_n |} \right )
\right ) - \log^2 x_n \right ) \right ) \nonumber \\
&=&\ \exp \left (-{1 \over 2C^2} \left (
\left ( \log x_n + \log \left ( 1 - {C^2\log k 
\over |\log x_n|} \right )
\right )^2 - \log^2 x_n \right ) \right )\nonumber \\
&<&\ \exp \left ( -{1 \over 2C^2} \left (  
\left ( \log x_n - {C^2\log k \over |\log x_n|}\right )^2 
- \log^2 x_n \right ) \right ) 
\nonumber \\
&=&\ \exp \left ( -{1 \over 2C^2} \left ( 2C^2\log k  + 
{C^4\log^2 k \over \log^2 x_n} \right ) \right ) \nonumber \\
&<&\ \exp \left ( - \log k \right )\ =\ {1 \over k}. \nonumber
\end{eqnarray}

Now, if we let $\rho_1 = x_{n+1}$, where $x_n$ and $x_{n+1}$ satisfy the 
conclusion of the lemma above, and so that $x_{n+1}$ lies in
$(0,\rho_0)$, then 
$$
\rho\ >\ \rho_1 \left (1 + {2C^2\log k \over |\log \rho_1|} \right ),
$$
implies $r(\rho) \geq r(x_n)$, which implies that
$$
{r(\rho_1) \over r(\rho)}\ \leq\ {r(x_{n+1}) \over r(x_n)}\ <\ {1 \over k}.
$$
\bigskip

\section{Proof of Proposition \ref{lambda_prop}} \label{lambda_section}
\bigskip

We note that 
\begin{eqnarray} \label{lambdaxyz_equation}
p^2 \mu_p(W)\ &=&\ \#\{x,y,z \in W\ :\ x+y \equiv 2z \pmod{p}\} \nonumber \\
&=&\ {1 \over p} \sum_{|a| < p/2} f_{W}\left ( {a \over p} 
\right )^2 f_{W}\left ( {-2a \over p} \right ) \nonumber \\
&=&\ \Sigma_1\ +\ \Sigma_2, 
\end{eqnarray}
where $\Sigma_1$ is the contribution of the terms when we sum with
$|a| < 5vq^{1-\rho_1 \epsilon^2}$, and $\Sigma_2$ is the contribution coming
from terms with $5vq^{1-\rho_1 \epsilon^2} \leq |a| < p/2$.  

We will now show that $\Sigma_1$ gives a good approximation to the 
sum in (\ref{lambdaxyz_equation}).  First, we will require the following
two results:

\begin{proposition} \label{approx_prop}
For 
\begin{equation} \label{u_range}
5vq^{-\rho_1 \epsilon^2}\ \leq\ |u|\ \leq\ {1 \over 2},
\end{equation}
for $\epsilon > 0$ sufficiently small, and $q$ sufficiently
large, we have that
$$
|f_{W}(u)|\ <\ 2 \pi \left ( {\epsilon \over 
\rho_1} \right )^{1/3} |S|. \nonumber
$$
\end{proposition}

\begin{lemma} \label{pq_switch}
Suppose $|2b| < Q$ and that $p = q/2 + \delta$, where $Q |\delta| < q/3$.  
Then, we have that 
$$
f_{W} \left ( {2b \over q} \right )\ =\ f_{W} \left ( {b \over p}
\right )\ +\ O \left ( {\delta Q |S| \over q}\right ).
$$
\end{lemma}
Note:  The proof of proposition \ref{approx_prop} can be found in section
\ref{approx_section}, and the proof of lemma \ref{pq_switch} can be 
found in section \ref{lemma_section}.  

We now establish that
\begin{equation} \label{2ap_bound}
5vq^{1-\rho_1 \epsilon^2} \leq |a| < p/2\ \Longrightarrow\ 
f_{W}\left ( {-2a \over p} \right )\ =\ O \left ( \left ( {\epsilon 
\over \rho_1} \right )^{1/3}|S| \right ).
\end{equation}
The proof goes as follows:  Suppose $5vq^{1-\rho_1 \epsilon^2} \leq |a| < p/2$, 
and let $b$ be the number in $(-p/2,p/2)$ that is congruent to 
$-2a$ modulo $p$.  If $|b| < 5vq^{1-\rho_1 \epsilon^2}$, then 
$b$ must be odd, and so from lemma \ref{pq_switch} we have that
$$
f_{W} \left ( {-2a \over p} \right )\ =\ 
f_{W} \left ( {b \over p} \right )\ =\ f_{W} \left ( {2b \over q} \right )
+ O \left ( q^{1 - \rho_1 \epsilon^2/4} \right )\ =\ O(\epsilon |S|),
$$ 
where the last inequality follows from (\ref{S''_equation}) since $b$
is odd.  If $|b| \geq 5vq^{1-\rho_1 \epsilon^2}$, 
then we have from proposition \ref{approx_prop} that 
$f_{W}(b/p)$ is $O( (\epsilon \rho_1^{-1})^{1/3} |S|)$.  Thus,
(\ref{2ap_bound}) follows.

From (\ref{2ap_bound}) and Parseval's identity we get
\begin{eqnarray} \label{sigma2_estimate}
|\Sigma_2|\ &=&\ O \left ( {|S| \over p} \left ( {\epsilon \over \rho_1} \right )^{1/3}
\right )\ \sum_{5q^{1-\rho_1 \epsilon^2} \leq |a| < p/2}
\left | f_{W} \left ( {a \over p} \right ) \right |^2 \nonumber \\
&=&\ O \left ( |S|^2 \left ( {\epsilon \over \rho_1} \right )^{1/3} \right ).
\end{eqnarray}

Applying corollary \ref{parseval_corollary}, 
which appears just after the statement of Parseval's identity in section 
\ref{lemma_section} and makes use of Parseval's identity and Cauchy's
inequality, together with lemma \ref{pq_switch},
we get that
\begin{eqnarray} \label{sigma1_estimate}
\Sigma_1\ &=&\ {1 \over p} \sum_{|a| \leq 5vq^{1-\rho_1 \epsilon^2}}
f_{W} \left ( {a \over p} \right )^2 
\left ( f_{W} \left ( {-4a \over q} \right ) + 
O(q^{1- \rho_1 \epsilon^2/4}) \right ) \nonumber \\
&=&\ {1 \over p} \sum_{|a| \leq 5vq^{1-\rho_1 \epsilon^2}}
f_{W} \left ( {a \over p} \right )^2 f_{S''} \left ( {-4a \over q}
\right )\ +\ O \left (q^{2 - \rho_1 \epsilon^2/4}\right ) \nonumber \\
&=&\ {1 \over p} \sum_{|a| \leq 5vq^{1-\rho_1 \epsilon^2}}
f_{W} \left ( {a \over p} \right ) f_{W} \left ( {-4a \over q} \right )
\nonumber \\
&&\ \ \ \ \ \times\ \left ( f_{W} \left ({2a \over q} \right ) 
\ +\ O(q^{1-\rho_1 \epsilon^2/4}) \right ) + O \left ( q^{2 - \rho_1 \epsilon^2/4} \right )
\nonumber \\
&=&\ {1 \over p} \sum_{|a| \leq 5vq^{1-\rho_1 \epsilon^2}}
f_{W} \left ( {a \over p} \right ) f_{W} \left ( {-4a \over q} \right )
f_{W} \left ( {2a \over q} \right ) \nonumber \\
&&\ \ \ \ \ \ +\ O \left ( q^{2 - \rho_1 \epsilon^2/4} \right ) \nonumber \\
&=&\ {1 \over p} \sum_{|a| \leq 5vq^{1-\rho_1 \epsilon^2}}
f_{W} \left ( {2a \over q} \right ) f_{W} \left ( {-4a \over q} \right )
\nonumber \\
&&\ \ \ \ \ \times\ \left ( f_{W} \left ( {2a \over q} \right )\ +\ 
O(q^{1-\rho_1 \epsilon^2/4}) \right )\ +\ O \left ( q^{2 - \rho_1 \epsilon^2/4} 
\right ) \nonumber \\
&=&\ {1 \over p} \sum_{|a| \leq 5vq^{1-\rho_1 \epsilon^2}}
f_{W} \left ( {2a \over q} \right )^2 f_{W} \left ( {-4a \over q} \right )
\ +\ O \left ( q^{2-\rho_1 \epsilon^2/4} \right ) \nonumber \\
&=&\ {1 \over p} \sum_{|a| \leq 10vq^{1-\rho_1 \epsilon^2} \atop a\ {\rm even}} 
f_{W} \left ({a \over q} \right )^2 f_{W} \left ({-2a \over q} \right )
\ +\ O(q^{2-\rho_1 \epsilon^2/4}).
\end{eqnarray}

Let $J$ be the set of all integers $a$ where either 
$|a| \leq 10vq^{1-\rho_1 \epsilon^2}$
and $a$ is odd, or where $10vq^{1-\rho_1 \epsilon^2} < |a| < q/2$.  Note that
this set $J$ is the set of all integers $a$ ``missing'' from the sum in 
the last line of (\ref{sigma1_estimate}).  For each $a \in J$,
let $|b| < q/2$ be congruent to $-2a \pmod{q}$.  We will show that 
for each such $a \in J$, 
\begin{equation} \label{f2a_equation}
f_W \left ( {-2a \over q} \right )\ =\ f_W \left ( {b \over q} \right )
\ =\ O(\epsilon |S|).
\end{equation}
To see this, we first consider the case where $10vq^{1-\rho_1 \epsilon^2} 
< |a| < q/2$.  For this case, either $b$ is odd, or else 
$|b| > 10vq^{1-\rho_1 \epsilon^2}$.  In either case, we deduce 
(\ref{f2a_equation}) from (\ref{S''_equation}).
For the case $|a| \leq 10vq^{1-\rho_1 \epsilon^2}$, $a$ odd, we also get from 
(\ref{S''_equation}) that (\ref{f2a_equation}) holds, because $-2a$ is
not divisible by $4$.

Now, by Parseval's identity we have that
\begin{eqnarray} \label{J_estimate}
{1 \over p} \sum_{a \in J} f_{W}\left ( {a \over q} \right )^2 
f_{W} \left ( {-2a \over q} \right )\ &=&\ O \left ( {\epsilon |S| \over p} 
\sum_{|a| < q/2} \left | f_{W} \left ( {a \over q} \right )\right |^2 \right )
\nonumber \\
&=&\ O ( \epsilon |S|^2 ). 
\end{eqnarray}
Thus, from (\ref{sigma1_estimate}) and (\ref{J_estimate}) we get
\begin{eqnarray}
\Sigma_1\ &=&\ {1 \over p} \sum_{|a| < q/2} f_W \left ( {a \over q} \right )^2
f_W \left ( {-2a \over q} \right )\ +\ O \left ( \epsilon |S|^2 \right ) 
\nonumber \\
&=&\ {q^3 \over p} \mu_{q}(W)\ +\ O \left ( \epsilon |S|^2 \right ). \nonumber
\end{eqnarray} 
Combining this with our estimate for $\Sigma_2$ above, as well as 
(\ref{lambdaxyz_equation}), we have
$$
p^2 \mu_p(W)\ =\ \Sigma_1 + \Sigma_2\ =\ {q^3 \over p} \mu_q(W)\ +\ 
O \left (  \left ( {\epsilon \over \rho_1} \right )^{1/3} |S|^2 \right ),
$$
which proves the proposition.

\bigskip

\section{Proof of Proposition \ref{approx_prop}} \label{approx_section}
\bigskip

Suppose that $u$ satisfies (\ref{u_range}), and let $a$ be any integer
so that
\begin{equation} \label{au_inequality}
\left | u - {a \over q} \right | \leq {1 \over 2q}.
\end{equation}
Since the set $W$ satisfies (\ref{S''_equation}), we have that
\begin{equation} \label{ab_bound}
{\rm if\ } b \in {\mathbb Z},\ |b| < vq^{1-\rho_1 \epsilon^2}-1,\ {\rm then\ }
\left | f_{W} \left ( {a-b \over q} \right ) \right | \leq \epsilon |S|.
\end{equation}
One basic consequence of this fact is the following lemma, which is
proved in section \ref{lemma_section}:

\begin{lemma} \label{short_sum_lemma}
If $N$ and $H$ are non-negative integers such that 
$[N+1,N+H] \subseteq [0,q-1]$, $a$ satisfies (\ref{au_inequality}),
$\epsilon > 0$ is sufficiently small, and $q$ is sufficiently large
in terms of $\rho$ and $\epsilon$,
then we have that 
$$
\left | \sum_{s \in W \atop s \in [N+1,N+H]} 
e\left ( {s a \over q} \right ) \right |\ <\ 
2 |S| \left ( {\epsilon H \over \rho_1 q} \right )^{1/3}.
$$
\end{lemma} 

To finish the proof of our proposition, we apply this lemma together
with partial summation:  Let $\delta = u - a/q$, and observe that 
$|\delta| \leq 1/(2q)$.  Let 
$$
h(x)\ =\ \sum_{s \in W \atop s \leq x} e \left ( {s a \over q} \right ).
$$
Then, we have
\begin{eqnarray}
|f_{W}(u)|\ &=&\ \left |f_{W}\left ( {a \over q} + \delta \right ) \right |
\ =\ \left | \int_0^q e(\delta x) d h(x) \right | \nonumber \\
&=&\ \left | e(\delta x) h(x) \biggr |_0^{q}\ -\ 2\pi i \delta \int_0^q 
e(\delta x) h(x) dx\right | \nonumber \\
&\leq&\ \left |f_{W}\left ( {a \over q} \right ) \right |\ 
\ +\ 2\pi \delta \int_0^q |h(x)| dx \nonumber \\
&\leq&\ \left |f_{W}\left ( {a \over q} \right ) \right |\ +\ 
4\pi \delta |S| \left ( {\epsilon \over \rho_1 q} \right )^{1/3} \int_0^q 
x^{1/3} dx \nonumber \\
&\leq&\ \left ( \epsilon\ +\ 3 \pi \delta q \left ( {\epsilon \over 
\rho_1} \right )^{1/3} \right ) |S|. \nonumber
\end{eqnarray} 
Using the fact that $|\delta| < 1/(2q)$, the proposition now follows.
\bigskip

\section{Technical Lemmas} \label{lemma_section}
\bigskip

In this section we will state a few technical lemmas that were
used throughout the paper, as well as provide proofs of these and 
other lemmas appearing in the paper. 
\bigskip

\begin{lemma} \label{poly_bound}
For $-1/2 \leq t \leq 1/2$, $t \neq 0$, we have
$$
\left | \sum_{j=N+1}^{N+H} e(jt) \right |\ \leq\ 
\min \left ( H, {1 \over 2|t|} \right ).
$$
\end{lemma}  

\begin{lemma} \label{parseval} {\bf (Parseval's Identity)}
If 
$$
f(t)\ =\ \sum_{j=0}^{q-1} \lambda_j e(jt),
$$
then
$$
\sum_{a=0}^{q-1} \left | f \left ( {a \over q} \right ) \right |^2
\ =\ q \sum_{j=0}^{q-1} |\lambda_j|^2.
$$
\end{lemma}
\bigskip

An almost immediate corollary of this lemma, which follows by combining
it with Cauchy's inequality, is as follows:

\begin{corollary} \label{parseval_corollary}
Suppose $W \subseteq \{0,1,...,q-1\}$, that $q/2 < p < 2q$,
and that both $b_1$ and $b_2$ are integers such that
$(b_1,q) = (b_2,p) = 1$.  Then, we have 
$$
\sum_{|a| < q/2} \left | f_W \left ( {b_1 a \over q} \right ) \right |
\ \left | f_W \left ( {b_2 a \over p} \right ) \right |\ =\ O(q|S|). 
$$
\end{corollary}
The proof of this result appears at the end of this section.
\bigskip

\noindent {\bf Proof of Corollary \ref{behrend_consequence}.}  
\bigskip

Set 
$$
L(x)\ =\ \exp( C\sqrt{\log x}),
$$
where $C$ is as given in theorem \ref{behrend_theorem}.  Let $x$ be the 
integer satisfying 
$$
4 L(x)\ <\ {1 \over \rho}\ \leq\ 4 L(x+1),
$$
and suppose that $q$ is any prime larger than $4x$.  Further, let
$S \subseteq \{1,2,...,x\}$ be any set of density at least $L(x)^{-1}$ 
having only trivial 3-term arithemtic progressions, as given by Theorem
\ref{behrend_theorem}.  

Define the set 
\begin{equation} \label{T_subset}
T \subseteq \left \{ 0,1,2,..., {q-1 \over 2} \right \} 
\end{equation}
as follows:
$$
T\ =\ \left \{s + 2kx\ :\ s \in S,\ 0 \leq k \leq K 
= \left \lfloor {q \over 4x} \right \rfloor \right \}.
$$
Note that
$$
{|T| \over q}\ =\ {|S| (K+1) \over q}\ >\ {|S| \over 4x}\ >\ 
{1 \over 4 L(x)}\ >\ \rho,
$$   
and so we see that $T$ contains density $> \rho$ of the residue classes
modulo $q$.  

We note that if $a,b,c \in T$, $0 \leq a,b,c \leq q-1$, then
$$
a+b = 2c\ \iff\ a+b \equiv 2c \pmod{q},
$$
since $T$ satisfies (\ref{T_subset}); also, since $S$ contains only
trivial 3-term arithmetic progressions, we have that  
$$
a+b = 2c\ \iff\ a = s + 2xk,\ b = s + 2x(k+2d),\ c = s + 2x(k+d),
$$
where $s \in S$.  
Thus, the number of triples $a,b,c \in T$ satisfying $a+b = 2c$ is
at most 
\begin{eqnarray}
|S|\ \#\{ k,d\ :\ 0 \leq k < k+d < k+2d < K\}\ &<&\ |S|K^2
\ \leq\ {q^2 \over x+1} \nonumber \\
&\leq&\ {q^2 \over \exp \left ( {1 \over C^2} \log^2 (4\rho) \right )},
\nonumber
\end{eqnarray}
which proves the corollary.
\bigskip

\noindent {\bf Proof of Lemma \ref{cube_lemma}.}
\bigskip

The proof is via the pigeonhole principle:  Let $a_1,...,a_t$ be all
the integers in $(0,q/2)$ such that 
\begin{equation} \label{size_condition}
\left | f \left ( {a \over q} \right ) \right |\ >\ \epsilon |S|,
\end{equation}
for $a=a_1,...,a_t$.  

We have that $t \leq (\rho_1 \epsilon^2)^{-1}/2$.  To see this,
first note that $|f_S(a/q)| = |f_S(-a/q)|$, and so the number of 
integers $a$ in $(0,q/2)$ satisfying (\ref{size_condition}) is 
at most half the total number of integers $a$ with $|a| < q/2$
satisfying (\ref{size_condition}), and this total number we know
to be at most $(\rho_1 \epsilon^2)^{-1}$.

Now, we note that to prove the lemma, it suffices
to find an integer $1 \leq j \leq q-1$ such that if $b_1,...,b_t$ are
the smallest numbers in absolute value that are congruent to 
$j a_1,...,j a_t$ modulo $q$, respectively, then 
\begin{equation} \label{bi_equation}
|b_i| \leq q^{1-\rho_1 \epsilon^2},\ {\rm for\ all\ }i=1,2,...,t.
\end{equation}
For if so, then if we let $h \equiv j^{-1} \pmod{q}$ and
$S' = h S$, then we get that
\begin{eqnarray}
&&\left |f_S \left ({hb \over q} \right ) \right | = 
\left |f_{S'}\left ( {b \over q} \right )\right | > 
\epsilon |S| = \epsilon |S'|\nonumber \\
&&\hskip0.5in \Longrightarrow\ hb \equiv \pm a_i \pmod{q}
\ {\rm where\ }i=1,2,...,t,\ \nonumber \\
&&\hskip1in {\rm or\ } b \equiv 0 \pmod{q}; \nonumber
\end{eqnarray}
but then this would mean that either $b \equiv 0 \pmod{q}$
or 
$$
b \equiv \pm h^{-1} a_i \equiv \pm j a_i \equiv \pm b_i \pmod{q},
$$
which means that the least residue in absolute value of $b$ modulo $q$ is 
$\leq q^{1-\rho_1 \epsilon^2}$.  

We now show how to find an integer $1 \leq j \leq q-1$ so that 
(\ref{bi_equation}) holds:  Partition the cube $[0,q-1]^t$
into the sub-cubes 
\begin{eqnarray}
&& [j_1 q^{1-\rho_1 \epsilon^2},\ (j_1+1)q^{1-\rho_1 \epsilon^2}]\ \times\  
[j_2 q^{1-\rho_1 \epsilon^2},\ (j_2 + 1)q^{1-\rho_1 \epsilon^2}] \nonumber \\
&& \hskip0.5in \times\ \cdots\ \times\ 
[j_t q^{1-\rho_1 \epsilon^2}, (j_t + 1)q^{1-\rho_1 \epsilon^2}], \nonumber
\end{eqnarray}
where $0 \leq j_1,...,j_t < q^{\rho_1 \epsilon^2}$.  Clearly, there are at most 
$$
\left ( q^{\rho_1 \epsilon^2} + 1 \right )^t\ <\ q
$$
such sub-cubes.  Now, consider the sequence
$$
(r a_1\ {\rm mod\ } q,\ ra_2\ {\rm mod\ }q,\ ...,\ ra_t\ {\rm mod\ }q),\ 
{\rm where\ } 0 \leq r \leq q-1.
$$
Since this sequence contains $q$ terms that lie inside the box $[0,q-1]^t$,
we must have that at least two of these terms lie in the same sub-cube.
If $r=r_1$ and $r=r_2$ are two such terms that correspond to points lying
in the same sub-cube, then it follows that
$$
( (r_1 - r_2) a_1\ {\rm mod\ q},\ ...,\ 
(r_1 - r_2) a_t\ {\rm mod\ q} )\ \in\ [-q^{1-\rho_1 \epsilon^2}, 
q^{1-\rho_1\epsilon^2}]^t,
$$
where here we take the least residue in absolute value for the entries.
So, letting $j = |r_1 - r_2|$ satisfies (\ref{bi_equation}).
We now have that (\ref{S'_property}) follows for this choice of $j$
(and $h$).
\bigskip
   
\noindent {\bf Proof of Lemma \ref{further_consequence}.}
\bigskip

If the conclusion of the lemma were false, then there exists $N \geq 1$
so that if $n > N$, then 
$$
{r(x_{n+1}) \over r(x_n)}\ \geq\ \exp \left ( -{1 \over 2C^2}
\left ( \log^2 x_{n+1} - \log^2 x_n \right ) \right ).
$$
Thus, for any $n > N$, 
\begin{equation} \label{xn_equation}
r(x_n)\ =\ r(x_N) \prod_{j=N+1}^n {r(x_j) \over r(x_{j-1})}
\ \geq r(x_N) \exp \left ( -{1 \over 2C^2} \left ( \log^2 x_n - 
\log^2 x_N \right ) \right ).  
\end{equation}
Noting here that $r(x_N) > 0$, which follows from (\ref{bourgain_varnavides})
from the introduction, we arrive at a 
contradiction, because from proposition \ref{behrend_consequence} we get that 
$$
r(x_n)\ <\ \exp \left ( -{1 \over C^2} \log^2 x_n \right ),
$$
which cannot be consistent with (\ref{xn_equation}) once $x_n$
is small enough (that is, once $n$ is large enough), 
Thus, the lemma follows.
\bigskip

\noindent {\bf Proof of Lemma \ref{pq_switch}.}
\bigskip

If $|2b| \leq Q$, and if $0 \leq s \leq q-1$,
then we have that
\begin{eqnarray}
e \left ( {2b s \over q} \right )
\ &=&\ e \left ( {2b s \over 2p - 2\delta} \right )
\ =\ e \left ( 2bs \left ( {1 \over 2p} + {2\delta \over 2p(2p - 2\delta)}
\right ) \right ) \nonumber \\
&=&\ e \left ( {b s \over p} \right ) 
\left ( 1\ +\ O \left ( {bs\delta \over p (2p - 2\delta)} \right ) \right ) \nonumber \\
&=&\ e \left ( {b s \over p} \right )\ +\ 
O \left ( {\delta Q \over q} \right ). \nonumber 
\end{eqnarray}
So,
\begin{eqnarray}
f_{W} \left ( {2b \over q} \right )\ &=&\ 
\sum_{s \in W} e \left ( {2b s \over q} \right )
\ =\ \sum_{s \in W} \left ( e \left ( {bs \over p} \right )\ +\ 
O \left ( {\delta Q \over q} \right ) \right ) \nonumber \\ 
&=&\ f_{W} \left ( {b \over p} \right )\ +\ O \left ( 
{\delta Q |S| \over q} \right ), \nonumber
\end{eqnarray}
which proves the lemma.
\bigskip

\noindent {\bf Proof of Lemma \ref{short_sum_lemma}.}
\bigskip

For $[N+1,N+H] \subseteq [0,q-1]$, let
$$
g(t)\ =\ \sum_{s \in W \atop s \in [N+1,N+H]} e(st),
$$
set
$$
D(t)\ =\ \sum_{j=N+1}^{N+H} e(jt),
$$
and let 
\begin{equation} \label{K_bound}
0 \leq K < vq^{1 - \rho_1 \epsilon^2} - 1 < q^{1- 3\rho_1\epsilon^2/4} - 1
\end{equation}
be some parameter, which is to be chosen later.

Then, we have for $a$ satisfying (\ref{au_inequality})
and $u$ satisfying (\ref{u_range}),

\begin{eqnarray} \label{g_value}
\left | g\left ( {a \over q} \right )\right |\ &=&\ 
\left | {1 \over q} \sum_{|b| < q/2} 
D\left ( {b \over q} \right ) f\left ({a-b \over q} \right ) \right | 
\nonumber \\
&\leq&\ {1 \over q} \sum_{|b| \leq K} \left | 
D \left ( {b \over q} \right ) \right |\ 
\left | f \left ( {a - b \over q} \right ) \right |
\ +\ \Sigma \nonumber \\
&\leq&\ {2 \epsilon K H |S| \over q} + \Sigma, 
\end{eqnarray}
where 
$$
\Sigma\ =\ {1 \over q} \sum_{K < |b| \leq q/2} \left | D \left ( {b \over q}
\right ) \right |\ \left | f \left ( {a-b \over q} \right ) \right |.
$$
We note that the last line of (\ref{g_value}) follows from 
(\ref{ab_bound}). 

To bound $\Sigma$ from above we will use Cauchy's inequality, Parseval's
identity, and the upper bound for $|D(t)|$ given by Lemma 
\ref{poly_bound} (which appears at the beginning of this section). 
We have
\begin{eqnarray}
\Sigma\ &\leq&\ {1 \over q} \left ( \sum_{K < |b| \leq q/2} 
\left | D \left ( {b \over q} \right ) \right |^2 \right )^{1/2} 
\left ( \sum_{K < b \leq q/2} \left | f \left ( {a-b \over q} \right )
\right |^2 \right )^{1/2} \nonumber \\
&<&\ {1 \over q} \left ( \sum_{|b| > K} {q^2 \over 4b^2} \right )^{1/2}
\left ( \sum_{j=0}^{q-1} \left | f \left ( {j \over q} \right ) \right |^2
\right )^{1/2} \nonumber \\
&<&\ {1 \over 2} \sqrt{q|S| \over K}
\ \leq\ {|S| \over 2} \sqrt{1 \over \rho_1 K}. \nonumber    
\end{eqnarray}

The value of $K$ that minimizes the last line of (\ref{g_value}) is
$$
K\ =\ \left ( {q^2 \over 16 \rho_1 \epsilon^2 H^2} \right )^{1/3},
$$
and we note that for $q$ sufficiently large and $\epsilon > 0$ sufficiently
small this will satisfy (\ref{K_bound}); and, with this choice of $K$, 
we get that
$$
\left | g \left ( {a \over q} \right ) \right |\ <\ 
2 |S| \left ( {\epsilon H \over \rho_1 q} \right )^{1/3},
$$
which proves the lemma.
\bigskip

\noindent {\bf Proof of Lemma \ref{poly_bound}.}  From the geometric
series identity, we have for $t \neq 0$, 
\begin{eqnarray}
\left | \sum_{j=N+1}^{N+H} e(jt) \right |\ &=&\ \left | {e(Ht) - 1 
\over e(t)-1} \right |\ \leq\ {2 \over |e(t/2) - e(-t/2)|} \nonumber \\
&=&\ {1 \over \sin(\pi |t|)}\ \leq\ {1 \over 2|t|}. \nonumber 
\end{eqnarray} 
The last inequality follows from the fact that for $0 \leq u \leq \pi/2$,
$$
\sin(u)\ \geq\ {2u \over \pi}.
$$
\bigskip

\noindent {\bf Proof of Corollary \ref{parseval_corollary}.}  
From Parseval's identity and Cauchy's inequality we have:
\begin{eqnarray} \label{parseval_cauchy}
&& \sum_{|a| < q/2} \left | f_W \left ( {b_1 a \over q} \right ) \right |
\ \left | f_W \left ( {b_2 a \over p} \right ) \right | \nonumber \\
&&\ \ \ \leq\ \left ( \sum_{|a| < q/2} \left | f_W \left ( {b_1 a \over q} \right ) \right |^2 
\right )^{1/2} \left ( \sum_{|a| < q/2} \left | f_W \left ( {b_2 a \over p} 
\right ) \right |^2 \right )^{1/2} \nonumber \\
&&\ \ \ \leq\ \left ( \sum_{|a| < q/2} \left | f_W \left ( {a \over q} \right ) \right |^2
\right )^{1/2} \left ( 2 \sum_{|a| < p/2} \left | f_W \left ( {a \over p} \right )
\right |^2 \right )^{1/2} \nonumber \\
&&\ \ \ \leq\ \left ( q |S| \right )^{1/2} \left ( 8p |S| \right )^{1/2}
\ =\ O(q|S|).
\end{eqnarray}
\bigskip

\section{Acknoledgements}
\bigskip

I would like to thank Ben Green for pointing out the reference
\cite{varnavides} below in an email many months ago, related to an
earlier paper of mine.

\end{document}